\documentclass[12pt]{article}
\usepackage{amssymb,amsmath}
\usepackage{amsthm}
\newcommand{\ff}{\mathbb{F}}
\newcommand{\ffq}{\mathbb{F}_q}
\newcommand{\nn}{\mathbb{N}}
\newcommand{\pp}{\mathbb{P}}
\newcommand{\zz}{\mathbb{Z}}
\newcommand{\qq}{\mathbb{Q}}
\newcommand{\rr}{\mathbb{R}}

\newcommand{\be}{\begin{enumerate}}
\newcommand{\ee}{\end{enumerate}}
\newcommand{\beq}{\begin{equation}}
\newcommand{\eeq}{\end{equation}}
\newcommand{\bea}{\begin{eqnarray}}
\newcommand{\eea}{\end{eqnarray}}
\newcommand{\beas}{\begin{eqnarray*}}
\newcommand{\eeas}{\end{eqnarray*}}
\newcommand{\st}{\,:\,}
\newcommand{\modd}[1]{\,(\mathrm{mod}\, #1)}
\newcommand{\bm}[1]{{\mbox{\boldmath $#1$}}}
\newcommand{\bms}[1]{{\mbox{\boldmath\scriptsize $#1$}}}
\newcommand{\bsp}{\!\!}

\def\({\left(}
\def\){\right)}

\swapnumbers
\theoremstyle{plain}
\newtheorem{theorem}{Theorem}[section]
\newtheorem{proposition}[theorem]{Proposition}
\newtheorem{lemma}[theorem]{Lemma}
\newtheorem{corollary}[theorem]{Corollary}

\theoremstyle{definition}
\newtheorem{defn}[theorem]{Definition}
\newtheorem{example}[theorem]{Example}

\theoremstyle{remark}

\begin{document}

\title{\bf Polynomial Coefficient Enumeration}
\author{Tewodros Amdeberhan\footnote{Department of Mathematics, Massachusetts Institute of Technology, Cambridge, MA 02139 Email: tewodros@math.mit.edu} 
\\  Richard P. Stanley\footnote{Department of Mathematics, Massachusetts Institute of Technology, Cambridge, MA 02139 Email: rstan@math.mit.edu}}
\date{November 21, 2008}

\maketitle
%\topmatter
%\begin{centering}
%\textbf{\large Polynomial Coefficient Enumeration}\\[.1em]
%(version of July 21, 2008)\\[1em]
%Tewodros Amdeberhan\\
%Richard P. Stanley\\[2em]
%\end{centering}
{\small \textsc{\bf Abstract}  Let $f(x_1,\dots,x_k)$ be a polynomial
  over a field $K$. This paper considers such questions as the enumeration
  of the number of nonzero coefficients of $f$ or of the number of
  coefficients equal to $\alpha\in K^*$. For instance, if $K=\ff_q$
  then a matrix formula is obtained for the number of coefficients of
  $f^n$ that are equal to $\alpha\in \ff_q^*$, as a function of
  $n$. Many additional results are obtained related to such areas as
  lattice path enumeration and the enumeration of integer points in
  convex polytopes. }

\vskip 5mm
%{\bf Key words:}
%{\bf MSC codes:}
\bigskip
%\document

\section{Introduction.} \label{sec1}

Given a polynomial $f\in\zz[x_1,\dots,x_n]$, how many coefficients
of $f$ are nonzero? For $\alpha\in\zz$ and $p$ prime, how many
coefficients are congruent to $\alpha$ modulo $p$? In this paper we
will investigate these and related questions. First let us review some
known results that will suggest various generalizations.

The archetypal result for understanding the coefficients of a
polynomial modulo $p$ is \emph{Lucas' theorem} \cite{6}: if $n, k$ are
positive integers with $p$-ary expansions $n=\sum_{i\geq 0}a_ip^i$ and
$k=\sum_{i\geq 0}b_ip^i$, then
 $$\binom{n}k\equiv \binom{a_0}{b_0}\binom{a_1}{b_1}\cdots \modd{p}.$$
Thus, for instance, it immediately follows that the number of odd
coefficients of the polynomial $(1+x)^n$ is equal to $2^{s(n)}$, where
$s(n)$ denotes the number of $1$'s in the $2$-ary (binary) expansion
of $n$. In particular, the number of odd coefficients of
$(1+x)^{2^m-1}$ is equal to $2^m$ (which can also be easily proved
without using Lucas' theorem).  In Section~\ref{sec2} we will vastly
generalize this result by determining the behavior of the number of
coefficients equal to $\alpha\in\ff_q$ of multivariable 
polynomials $f(\bm{x})^n$ as a
function of $n$ for any $f(\bm{x})\in\ff_q[\bm{x}]$.  A few related results
for special polynomials $f(\bm{x})$ will also be given.

In Sections~\ref{sec3}--\ref{sec5} we turn to other multivariate
cases. The following items give a sample of our results on sequences of 
polynomials not of the form $f(\bm{x})^n$.

%In Section~\ref{sec6}, we conclude with 
%certain remarks and conjectures related to generalized 
%Chebychev polynomials that appear in Section~\ref{sec5}. 
%Some known results akin to our work are the following.

\be \item Extensions of Lucas' theorem and related results. For
instance, in Section~\ref{sec3} we determine the number of coefficients
of the polynomial $\prod_{i=1}^n(1+x_i+x_{i+1})$ that are not
divisible by a given prime $p$.

\item It is easy to see, as pointed out by E. Deutsch \cite{5}, that
  the number of nonzero coefficients of the polynomial
  \beq \prod_{i=1}^n(x_1+x_2+\cdots +x_i) \label{eq:cat} \eeq  
  is the Catalan number $C_n=\frac{1}{n+1}\binom{2n}{n}$.
  This statement is generalized in Section~\ref{sec4}, where we
  interpret the number of nonzero coefficients of more general
  polynomials in terms of lattice path counting.

\item A number of examples are known where the nonzero coefficients
  are in bijection with the lattice points in a polytope.  The
  machinery of counting lattice points can be brought to bear. For
  example, this technique has been used to show that the number of
  nonzero coefficients of the polynomial $\prod_{1\leq i<j\leq
    n}(x_i+x_j)$ is equal to the number of forests on an $n$-vertex
  set (equivalent to the case $n=1$ of \cite[Exer. 4.32(a)]{18}). For
  a more complicated example of this nature, see \cite{17}. In
  Section~\ref{sec4} we develop a connection between the nonzero
  coefficients of polynomials of the form
  $\prod_{i=1}^k(x_1+x_2+\cdots +x^{\lambda_i})$ and a polytope
  studied by Pitman and Stanley \cite{14}.

\item The number of odd coefficients of the polynomial $\prod_{1\leq
i<j\leq n}(x_i+x_j)$ is equal to $n!$. This result is most easily seen
using the trick
 \begin{align}
 \prod_{1\leq i<j\leq n}(x_i+x_j)&\equiv \prod_{1\leq i<j\leq
 n}(x_i-x_j) \modd{2} \nonumber\\
 &=\sum_{w\in S_n}\text{sgn}(w)x_{w(1)}^{n-1}x_{w(2)}^{n-2}\cdots
 x_{w(n-1)}, \nonumber\end{align}
by the expansion of the Vandermonde determinant. An alternative
cancellation proof was given by Gessel \cite{9}, i.e., pairing off equal
monomials among $2^{\binom{n}2}$ monomials appearing in the product
until $n!$ distinct monomials remain. Some other results along this
line appear in Section~\ref{sec5}.
 \ee

\section{Powers of polynomials over finite fields}
\label{sec2}The prototype for the next result is the fact mentioned in
Section~\ref{sec1} that the number of odd coefficients of the
polynomial $(1+x)^{2^n-1}$ is equal to $2^n$.

Let $q=p^r$ where $p$ is a prime and $r\geq 1$. Let $k\geq 1$ and set
$\bm{x}=(x_1,\dots,x_k)$. Fix $f(\bm{x})\in\ff_q[\bm{x}]$ and
$\alpha\in\ff_q^{*}$. Define $N_\alpha(n)$ to be the number of
coefficients of the polynomial $f(\bm{x})^n$ that are equal to $\alpha$.
Let $\nn=\{0,1,2,\dots\}$. We use the multivariate notation
$\bm{x}^\gamma=x_1^{\gamma_1}\cdots x_k^{\gamma_k}$, where
$\gamma=(\gamma_1,\dots,\gamma_k)\in \nn^k$.

\begin{theorem} \label{thm:fxkn} There exist $\zz$-matrices $\Phi_0,
  \Phi_1,\dots,\Phi_{q-1}$ of some square size, and there exist a row
  vector $u$ and a column vector $v$ with the following property. For
  any integer $n\geq 1$ let $a_0+a_1 q+\cdots +a_r q^r$ be its base
  $q$ expansion, so $n=a_0+a_1 q+\cdots +a_r q^r$ and $0\leq a_i\leq
  q-1$.  Then
   $$ N_\alpha(n) = u\Phi_{a_r}\Phi_{a_{r-1}}\cdots \Phi_{a_0}v. $$
The vector $v$ and matrices $\Phi_i$ do not depend on $\alpha$.
\end{theorem}

\proof
Our proof is an adaptation of an argument of Moshe
\cite[Thm.~1]{moshe}. Suppose that $a_0,a_1,\dots$ is an infinite
sequence of integers satisfying $0\leq a_i\leq q-1$. Let ${\cal P}'$
be the Newton polytope of $f$, i.e.,
the convex hull in $\rr^k$ of the exponent vectors of monomials
appearing in $f$, and let ${\cal P}$ be the convex hull of
${\cal P}'$ and the origin. If $c>0$, then write
$c{\cal P}=\{ cv\st v\in{\cal P}\}$. Set $S=(q-1){\cal P}\cap \nn^k$
and $r_m=\sum_{i=0}^m a_iq^i$.
%Set $[0,j]=\{0,1,\dots,j\}$. Choose j large enough so that
%${\cal P}\cap\nn^k\subseteq [0,j]^k:=S$. Thus if $c\in\pp$ then $c{\cal
%  P}\subseteq cS=[0,cj]^k$.

Suppose that $f(\bm{x})^{r_m} =\sum_\gamma c_{m,\gamma}
\bm{x}^\gamma$. We set $f(\bm{x})^{r_{-1}}=1$. Let $\ffq^S$ be the set
of all functions $F\colon S\rightarrow \ffq$.  We will index our
matrices and vectors by elements of $\ffq^S$ (in some order).
%and
%let ${\cal V}$ be the real vector space of all functions
%$\varphi\colon \ffq^S \rightarrow \rr$.
Set
  $$ R_m= \{0,1,\dots,q^{m+1}-1\}^k. $$
For $m\geq -1$, define a column vector $\psi_m$ by letting $\psi_m(F)$
(the coordinate of $\psi_m$ indexed by $F\in\ffq^S$) be the number of
vectors $\gamma\in R_m$ such that for all $\delta\in S$ we have
$c_{m,\gamma+q^{m+1}\delta}=F(\delta)$.
%In particular,
%    $$ \psi_{-1}(F) = \left\{ \begin{array}{rl}
%        1, & \mathrm{if}\ F(\beta) = \delta_{\beta, (0,0,\dots,0)}\\
%        0, & \mathrm{otherwise}. \end{array} \right. $$
Note that by the definition of $S$ we have
$c_{m,\gamma+q^{m+1}\delta}=0$ if $\delta\not\in S$. (This is the
crucial finiteness condition that allows our matrices and vectors to
have a fixed finite size.) Note also that given $m$, every point $\eta$
in $\nn^k$ can be written uniquely as $\eta=\gamma+q^{m+1}\delta$ for
$\gamma\in R_{m+1}$ and $\delta$ in $\nn^k$.

For $0\leq i\leq q-1$ define a matrix $\Phi_i$ with rows and columns
indexed by $\ffq^S$ as follows. Let $F,G\in\ffq^S$. Set
   \beas g(\bm{x}) & = & f(\bm{x})^i\sum_{\beta\in S}
   G(\beta)\bm{x}^\beta\\
  & = & \sum_\gamma d_\gamma \bm{x}^\gamma \in\ffq[\bm{x}]. \eeas
Define the $(F,G)$-entry $(\Phi_i)_{FG}$ of $\Phi_i$ to be the number
of vectors $\gamma\in R_0=\{0,1,\dots,q-1\}^k$ such that for all
$\delta\in S$ we have $d_{\gamma+q\delta}=F(\delta)$.
%[check??]
A straightforward computation shows that
   \beq \Phi_{a_m}\psi_{m-1} = \psi_m,\ \ m\geq 0. \label{eq:phipsi}
   \eeq
Let $u_\alpha$ be the row vector for which $u_{\alpha}(F)$ is the number of
values of $F$ equal to $\alpha$, and let $n=a_0 +a_1q+\cdots +a_rq^r$
as in the statement of the theorem. Then it follows from
equation~\eqref{eq:phipsi}
that
  $$ N_\alpha(n) = u\Phi_{a_r}\Phi_{a_{r-1}}\cdots
  \Phi_{a_0}\psi_{-1}, $$
completing the proof. \qed

%Now $f(\bm{x})^{r_m}=f(\bm{x})^{r_{m-1}}f(\bm{x}^{q^m})^{a_m}$,
%where $\bm{x}^{q^m} =
% (x_1^{q^m},\cdots,x_k^{q^m})$.  From this formula

% Let us
% illustrate this contention with an example which should make the
% general argument clear.
%
% Let $q=3$ and $f=x+2y+x^2y$.

% Indeed, for any $g\in\ffq[x]$ and any $\gamma\in R_m$, let
% $g_{m,\gamma}\in \ffq[x]$ consist of those terms $d_\beta x^\beta$ of
% $g$ for which $\beta= \gamma+q^{m-1}\delta$ for some $\delta\in\nn^k$.

% Now let $n=a_0+a_1q+\cdots+a_rq^r$ as in the statement of the
% theorem. Note that $N_\alpha(n)$ can be obtained from $\Phi_{}$

\begin{example}
We illustrate the above proof with the simplest possible example,
since any more complicated example involves much larger matrices. Take
$q=2$, $k=1$, $f(x)=x+1$ and (necessarily) $\alpha=1$. Then
$S=\{0,1\}$. A function $F\colon S\rightarrow \mathbb{F}_2$
will be identified with the binary word $f(0)f(1)$, and our vectors
and matrices will be indexed by the words 00, 10, 01, 11 in that
order. Take $a_0 = 1$, $a_1=1$, $a_2=0$, $a_3=1$. Then
  \beas (1+x)^0 & = & 1\\
     (1+x)^1 & = & 1+ x\\
     (1+x)^{1+2^1} & = & 1+x+x^2+x^3\\
     (1+x)^{1+2^1+0\cdot 2^2} & = & 1+x+x^2+x^3\\
     (1+x)^{1+2^1+0\cdot 2^2+2^3} & = & 1+x+x^2+x^3+x^8+x^9+x^{10}
                         +x^{11} \eeas
%     (1+x)^{1+2^1+0\cdot 2^2+2^3+0\cdot 2^4} & = &
%     1+x+x^2+x^3+x^8+x^9+x^{10} +x^{11} \eeas

  $$ \psi_{-1} = \left[ \begin{array}{c} 0\\ 1\\  0\\ 0
      \end{array} \right],\
     \psi_0 = \left[ \begin{array}{c} 0\\ 2\\ 0\\ 0
      \end{array} \right],\
     \psi_1 = \left[ \begin{array}{c} 0\\ 4\\ 0\\ 0
      \end{array} \right],\
     \psi_2 = \left[ \begin{array}{c} 4\\ 4\\ 0\\ 0
      \end{array} \right],\
     \psi_3 = \left[ \begin{array}{c} 8\\ 8\\ 0\\ 0
      \end{array} \right] $$

   $$ \Phi_0=\left[ \begin{array}{cccc}
     2 & 1 & 1 & 0\\ 0 & 1 & 1 & 2\\ 0 & 0 & 0 & 0\\ 0 & 0 & 0 & 0
     \end{array} \right], \ \
   \Phi_1=\left[ \begin{array}{cccc}
     2 & 0 & 0 & 1\\ 0 & 2 & 1 & 0\\ 0 & 0 & 1 & 0\\ 0 & 0 & 0 & 1
   \end{array} \right] $$
For instance, to obtain $\psi_2$, break up $(1+x)^{1+2^1+0\cdot 2^2}$
into $2^{2+1}=8$ parts according to the congruence class of the
exponent modulo $2^{2+1}=8$: $1+0\cdot x^8$, $x+0\cdot x^9$,
$x^2+0\cdot x^{10}$, $x^3+0\cdot x^{11}$, $0\cdot x^{4} +0\cdot
x^{12}$, $0\cdot x^5 +0\cdot x^{13} $, $0\cdot x^6 +0\cdot x^{14}$,
$0\cdot x^7 +0\cdot x^{15}$. Four of these ``sections'' have
coefficient sequence 00 and four have 10, so
$\psi_2(00)=\psi_2(10)=4$, $\psi_2(01)=\psi_2(11)=0$. To get the
$01$-column of $\Phi_1$ (or the third column using the usual indexing
$1,2,3,4$ of the rows and columns), multiply $x$ (corresponding to
$01$) by $f(x)=1+x$ to get $x+x^2$. Bisect $x+x^2$ into the two
sections $0+x^2$ and $x+0\cdot x^3$. The coefficient sequences of
these two sections are 01 and 10. Hence $\Phi_1(01)=\Phi_1(10)=1$,
$\Phi_1(00) = \Phi_1(11)=0$, and $u=[0,1,1,2]$.
\end{example}

\begin{corollary} \label{cor:mainmod}
Preserve the notation of Theorem~\ref{thm:fxkn}, and set
  $$ F_{f,\alpha}(t) =\sum_{m\geq 0}N_\alpha(1+q+\cdots+q^{m-1})t^m. $$
Then $F_{f,\alpha}(t)$ is a rational function of $t$.
\end{corollary}

\proof
By Theorem~\ref{thm:fxkn} we have
  $$ N_\alpha(1+q+\cdots+q^{m-1}) = u\Phi_1^mv. $$
The proof now follows from standard arguments (e.g.,
\cite[Thm.~4.7.2]{18}) from linear algebra.
\qed

\medskip
\textsc{Note.} It is clear from the proof of Theorem~\ref{thm:fxkn}
that Corollary~\ref{cor:mainmod} can be considerably generalized.  For
instance, for any
$f(\bm{x}),g(\bm{x})\in\ffq[\bm{x}]=\ffq[x_1,\dots,x_k]$, $r\geq 1$,
and $\alpha\in\ffq^*$, let $L(m)$ be the number of coefficients of the
polynomial $g(\bm{x})f(\bm{x})^{(q^{rm}-1)/(q^r-1)}$ equal to
$\alpha$. Then $\sum_{m\geq 0} L(m)t^m$ is a rational function of $t$.

The examples provided below (working in $\ff_2[\bm{x}]$) demonstrate the 
conclusion promised by Corollary~\ref{cor:mainmod}. In some cases we give
an independent argument to arrive at the generating function.

\begin{example} Suppose $\bm{x}=(x_1,\dots,x_k)$ and let 
$$f_n(\bm{x})=\left(1+\sum_{i=1}^kx_i+x_1
\sum_{i=2}^k x_i^2\right)^{2^n-1}.$$
Then the number of odd coefficients is
 $$ N_1(f_n(\bm{x}))=k\cdot(k+1)^n-(k-1)\cdot k^n.$$

 %[\textbf{Tewodros.} The proof below is hard to understand. You use
 %terms like row, NE, west-coastal, etc., without defining the array
 %they are applying to. And it was hard for me to understand the part
 %beginning with ``However, there are additional interactions between
 %the first and third rows, via the second.'' What exactly are the
 %interactions? How exactly does a certain part ``act as a catalyst''?] 

\proof We work in $\ff_2[\bm{x}]$. In 
particular, using $(a+b)^2=a^2+b^2$, it is evident that 
\begin{align*}
f_n(\bm{x})
&=((1+x_2+\cdots+x_k)+x_1(1+x_2^2+\cdots+x_k^2))^{2^n-1}\\
&=\sum_{j=0}^{2^n-1}x_1^j\binom{2^n-1}{j}(1+x_2+\cdots+x_k)^{2^n-1-j}
(1+x_2^2+\cdots+x_k^2)^j\\
&=\sum_{j=0}^{2^n-1}x_1^j(1+x_2+\cdots+x_k)^{2^n-1-j}
(1+x_2^2+\cdots+x_k^2)^j\\
&=\sum_{j=0}^{2^n-1}x_1^j(1+x_2+\cdots+x_k)^{2^n-1+j}.
\end{align*}
Therefore our enumeration translates to
\beq
N(f_n(\bm{x}))=\sum_{j=0}^{2^n-1}N((1+x_2+\cdots+x_k)^{2^n-1+j}).
\label{eq:eqfull1} \eeq
We prove the assertion of this example by inducting on $n$. The
base case $n=0$ is obvious. Assume its validity for $n-1$. 

An argument verifying $N((1+x)^m)=2^{s(m)}$ (see
Introduction) also shows that $N((1+x_2+\cdots+x_k)^m)=k^{s(m)}$, where
$s(m)$ denotes the number of $1$'s in the binary expansion of $m$.
Thus equation~\eqref{eq:eqfull1} becomes
\beq
N(f_n(\bm{x}))=\sum_{j=0}^{2^n-1}k^{s(2^n-1+j)}.\label{eq:eqfull2} \eeq
Applying the following simple facts
$$s(2^n-1+j)=\begin{cases}n,& j=0\\
s(2^{n-1}-1+j),& 1\leq j\leq 2^{n-1}-1\\
1+s(2^{n-1}-1+j^{\prime}),& 0\leq j^{\prime}:=j-2^{n-1}\leq 2^{n-1}-1
\end{cases} $$
in equation~\eqref{eq:eqfull2} produces
\begin{align*}
N(f_n(\bm{x}))&=\sum_{j=0}^{2^{n-1}-1}k^{s(2^n-1+j)}+
\sum_{j=2^{n-1}}^{2^n-1}k^{s(2^n-1+j)}\\
&=k^n+
\sum_{j=1}^{2^{n-1}-1}k^{s(2^{n-1}-1+j)}+
\sum_{j=0}^{2^{n-1}-1}k^{1+s(2^{n-1}-1+j)}\\
&=k^n-k^{n-1}+N(f_{n-1}(\bm{x}))+kN(f_{n-1}(\bm{x})).
\end{align*}
Then by the induction assumption, we get
$$N(f_n(\bm{x}))=k^n-k^{n-1}+
(k+1)(k(k+1)^{n-1}-(k-1)k^{n-1})=k(k+1)^n-(k-1)k^n.$$

\end{example}  \qed

\begin{example} Denote the forward shift by $E(a_n)=a_{n+1}$.
\be
\item[(i)] If $k=2$ and $f_n=(1+x_1+x_2+x_1x_2^2)^{2^n-1}$
%\in\ff_2[x_1,x_2]$ 
then
$$N(f_n)=2\cdot3^n-2^n.$$
The same formula holds for $g_n=(1+x_1+x_2^2+x_1x_2)^{2^n-1}$.

\item[(ii)] When $k=3$ so
  $f_n=(1+x_1+x_2+x_3+x_1x_2^2+x_1x_3^2)^{2^n-1}$, then $a_n=N(f_n)$  
satisfies 
 $$ (E^2-7E+12)a_n=(E-4)(E-3)a_n=0.$$ 
We get $a_n=3\cdot4^n-2\cdot3^n$. 

\item[(iii)] If $g_n=(1+x_1+x_2+x_3+x_1x_2^2+x_2x_3^2)^{2^n-1}$ then
  $a_n=N(g_n)$ satisfies  
 $$ (E^2-7E+10)a_n=(E-5)(E-2)a_n=0.$$ 
We get $a_n=\frac13(4\cdot5^n-2^n)$.

\item[(iv)] Let $g_n=(1+x_1+x_2+x_3+x_1x_2^2)^{2^n-1}$. Then
  $a_n=N(g_n)$ satisfies  
$(E^2-6E+7)a_n=0$. Here the eigenvalues are not even rational; i.e.
if $c=3+\sqrt{2}$ and $b=1+\sqrt{2}$ and we write $\overline{\mu}$ for
quadratic conjugation in  
$\zz[\sqrt{2}]$, then
$$a_n=\frac{bc^n+\overline{b}\overline{c}^n}2.$$

\item[(v)]Let $g_n=(1+x_1+x_2+x_3+x_1x_2+x_1x_3^2)^{2^n-1}$. Now
  $a_n=N(g_n)$ satisfies $(E^2-7E+8)a_n=0$.  
Again the eigenvalues are not even rational; i.e.
if $c=\frac{7+\sqrt{17}}2$ and $b=17+5\sqrt{17}$ and we write
$\overline{\mu}$ for quadratic  
conjugation in $\qq[\sqrt{17}]$, then
$$a_n=\frac{bc^n+\overline{b}\overline{c}^n}{34}.$$

\item[(vi)] Let $g_n=(1+x_1+x_2^2+x_1x_2^3)^{2^n-1}$, $c=2+\sqrt{3}$,
  and $\overline{c}=2-\sqrt{3}$. We have
 $$ N(g_n)=\frac{c^{n+1}+\overline{c}^{n+1}-2^n}3.$$

\item[(vii)] Let $g_n=(1+x_1+x_2+x_1^2x_2^2)^{2^n-1}$. Then
  $a_n=N(g_n)$ satisfies  
$(E^4-5E^3+6E^2-2E-4)a_n=0$ whose eigenvalues are complex. 
Put $c=\frac{3+\sqrt{17}}2$ and $b=17^2+73\sqrt{17}$ with conjugates in 
$\qq[\sqrt{17}]$. Also let $i=\sqrt{-1}$, then
$$N(g_n)=\frac{bc^n+\overline{b}\overline{c}^n}{442}-
\frac1{13}(-2+3i)(1+i)^n-\frac1{13}(2+3i)(1-i)^n.$$
%\frac{2}{\sqrt{13}}2^{n/2}\cos\left(\frac{\pi n}4-\arctan(3/2)\right).$$

%[\textbf{Tewodros:} I think that the above formula would be easier to
%understand if expressed in terms of $\frac 12(3\pm\sqrt{17})$ and $1\pm
%i$, even though it is no longer in ``real'' form. ]

\item[(viii)] Let $g_n=(1+x_1^2+x_2^2+x_1x_2^3)^{2^n-1}$. Then 
$$\sum_{n\geq 0}N(g_n)z^n=\frac{1-2z+4z^2}{1-6z+12z^2-12z^3}.$$

\item[(ix)] Let $g_n=(1+x_1+x_2+x_2^2)^{2^n-1}$. Then 
$$\sum_{n\geq 0}N(g_n)z^n=\frac{1+2z}{1-2z-4z^2}.$$
In fact $N(g_n)=2^nF(n+2)$ where $F(n)$ is the Fibonacci sequence;
that is, $F(1)=1, F(2)=1$, $F(m)=F(m-1)+F(m-2)$ for $m\geq 3$.

\medskip
\item[(x)] The following three statements are immediate from example (ix).

\item[(a)] Given $k\in\pp$, break up the binary digits of 
$k$ into maximal strings of consecutive $1$'s and let the lengths of these
strings be the multiset $\bm{k}=\{k_1,k_2,\dots\}$. Then we obtain the 
averaging value
$$\frac1{2^n}\sum_{k=0}^{2^n-1}\prod_{\bm{k}}\frac{2^{k_i+2}+(-1)^{k_i+1}}3=F(n+2).$$ 

\item[(b)] The probability of not landing two consecutive heads in a
  fair toss of $n$ coins is equal to that of finding
a 1 in the triangle of coefficients formed by 
$$(1+x_2+x_2^2)^m \mod 2,\qquad 0\leq m\leq 2^n-1,$$ 
which in fact is $F(n+2)/2^n$. 

\item[(c)] The generating function 
$$\Lambda(z)=\sum_{m=0}^{\infty}N((1+x_2+x_2^2)^m)z^m$$
satisfies $\Lambda(z)=(1+2z)\Lambda(z^2)$. 
\ee \end{example}

\begin{example} 
\be

\item[(i)] Let $g_n=(1+x_1+x_2+x_2^3)^{2^n-1}$. Then 
$$\sum_{n\geq 0}N(g_n)z^n=\frac{1+z-2z^3}{1-3z-2z^2+2z^3+4z^4}.$$

\item[(ii)] Let $g_n=(1+x_1+x_2+x_2^4)^{2^n-1}$. Then  
$$\sum_{n\geq 0}N(g_n)z^n=\frac{1+z+4z^2+2z^3-4z^4}{1-3z-2z^3-8z^4+8z^5}.$$

\item[(iii)] If $g_n=(1+x_1+x_2+x_3+x_1x_2^2+x_1x_3^2+x_2x_3^2)^{2^n-1}$, then 
$$\sum_{n\geq 0}N(g_n)z^n=\frac{(1-z)^2}{1-9z+23z^2-19z^3}.$$

\item[(iv)] If $g_n=(1+x_1+x_2+x_3+x_1x_2^2+x_2x_1^2)^{2^n-1}$, then
$$\sum_{n\geq 0}N(g_n)z^n=\frac{1-z+2z^2-4z^3}{1-7z+12z^2-12z^3+8z^4}.$$
\ee \end{example}

\begin{example}  

\emph{Symmetric polynomials}. Consider the Vandermonde-type 
polynomials in $\ff_2[\bm{x}]$ given by
$$V(k,n):=\prod_{1\leq i<j\leq k}(x_i-x_j)^{2^n-1} \qquad\text{and}$$
$$V^{\prime}(k,n):=\left(1+\prod_{1\leq i<j\leq k}(x_i-x_j)\right)^{2^n-1}.$$
We find that 

\item[(i)] $N(V(2,n))=2^n$ and $N(V^{\prime}(2,n))=3^n$.

\item[(ii)] If $c=\frac{5+\sqrt{33}}2, b=11+3\sqrt{33}$ with conjugation in
$\qq[\sqrt{33}]$, then 
$$N(V(3,n))=6\cdot4^{n-1}, \qquad N(V^{\prime}(3,n))=\frac{bc^n+\overline{b}\overline{c}^n}{22}.$$

\item[(iii)] $N(V(4,n))=5\cdot8^n-8\cdot2^n$. 

\end{example}

\medskip
Consider the special case of Corollary~\ref{cor:mainmod} that $k=1$ and
$f(x)$ has the form $g(x)^{q-1}$. Thus
$f(x)^{(q^n-1)/(q-1)} = g(x)^{q^n-1}$. For $\alpha\in\ffq^*$ let
$M_\alpha(m)$ be the number of coefficients of $g(x)^{q^m-1}$ equal to
$\alpha$. We can give more precise
information about the linear recurrence satisfied by $M_\alpha(m)$.
To give a slightly more general result, we also fix $c\in\pp$.
Without loss of generality we may assume $g(0)\neq 0$. For
$m\in\pp$ such that $q^m\geq c$, let $N_{\alpha}(m)$
denote the number of coefficients of the polynomial $g(x)^{q^m-c}$
that are equal to $\alpha$.

\begin{theorem} \label{thm:1}
\be\item[(a)] There exist periodic functions $u(m)$
and $v(m)$ depending on $g(x), c$, and $\alpha$, such that
  \beq N_{\alpha}(m)=u(n)q^m+v(m) \label{eq1} \eeq
for $m$ sufficiently large.
 \item[(b)] Let $d$ be the least positive integer for which $g(x)$
divides $x^{q^m(q^d-1)}-1$ for some $m\geq 0$. In other words, $d$ is
the degree of the extension field of $\ff_q$ obtained by adjoining
all zeros of $g(x)$. Then the functions $u(m)$ and $v(m)$ have period
$d$ (and possibly smaller periods, necessarily dividing $d$).
 \item[(c)] Let $\mu$ be the largest multiplicity of any
irreducible factor (or any zero) of $g(x)$. Then equation~\eqref{eq1}
holds  for all $m\geq \lceil{\log_q\mu c}\rceil$. In particular, if
$g(x)$ is
squarefree and $c=1$, then \eqref{eq1} holds for all $m\geq 0$.
 \item[(d)] If $g(x)$ is primitive over $\ff_q$ and $c=1$, then
$d=\deg f$ and $u(m)=dq^{d-1}/(q^d-1)$, a constant.
  \ee
\end{theorem}

\proof We have $g(x)^{q^m-c}=g(x^{q^m})/g(x)^c$. Let
$g(x)=a_0+a_1x+\cdots+a_{\delta}x^{\delta}$, and for $0\leq i<\delta$
set
 \begin{align} G_{in}(x)&=\frac{(a_0+a_1x+\cdots+a_ix^i)^{q^m}}{g(x)^c}\\
&=\frac{a_0+a_1x^{q^m}+a_2x^{2q^m}+\cdots+a_ix^{iq^m}}{g(x)^c}\\
&=g_{im}(x)+\frac{h_{im}(x)}{g(x)^c},\end{align}
where $g_{im}(x), h_{im}(x)\in\ff_q[x]$ and $\deg h_{im}(x) <
c\deg g(x)$. Thus $h_{im}(x)$ is the remainder upon dividing
$(a_0+a_1x+\cdots+a_ix^i)^{q^m}$ by $g(x)^c$. Hence $h_{im}(x)$
determines the coefficient of $x^j$ in $g(x)^{q^m-c}$ for $iq^m\leq
j<(i+1)q^m$. These coefficients will be periodic, with period of the
form $\pi=q^s(q^d-1)$ for some $s\geq 0$. If $\alpha$ occurs $k$ times
within each period, then the number of times $\alpha$ occurs as a
coefficient of $x^j$ in $g(x)^{q^m-c}$ for $iq^m\leq j<(i+1)q^m$ has
the form $kq^m/\pi+v_i(m)$, where $v_i(m)$ depends only on $q^m$
modulo $\pi$.

Suppose that $0\leq l<m$ and $l\equiv m\modd{d}$. Then
 $$ G_{im}(x)-G_{il}(x)=\frac{a_1(x^{q^m}-x^{q^l})+\cdots+
   a_i(x^{iq^m}-x^{iq^l})}{g(x)^c}\in\ff_q[x]. $$
Hence if $l$ is large enough so that $g(x)^c$ divides
$x^{q^l(q^d-1)}$, then it follows that $h_{il}(x)=h_{im}(x)$. Thus the
polynomial $h_{im}(x)$ depends only on the congruence class of $m$
modulo $d$ for $m\geq l$. We can take $l$ to be the least integer such
that $g(x)^c$ divides $x^{q^l(q^d-1)}$. Thus $l$ is the least integer
for which $\mu c\leq q^l$, i.e.,  $l=\lceil{\log_q\mu c}\rceil$.

From the above discussion it follows that the coefficients of
$g(x)^{q^m-c}$ are periodic between $1$ and $x^{q^m-1}$ (i.e., for the
coefficients of $x^0=1,x,x^2,\dots,x^{q^m-1}$), then periodic between
$x^{q^m}$ and $x^{2q^m-1}$, etc. The lengths of these periods can all
be taken to be $d$. (Of course $d$ may not be the length of the
\emph{minimal} period).  Moreover, the coefficients themselves within
each period depend only on $m$ modulo $d$. If the number of times
$\alpha\in\ff_q$ appears within each period between $x^{iq^m}$ and
$x^{(i+1)q^m-1}$ is $k$, then the total number of coefficients between
$x^{iq^m}$ and $x^{(i+1)q^m-1}$ that are equal to $\alpha$ is $kq^m/d$
plus an error that is periodic with period $d$. It follows that
$N_{\alpha}(m)=u(m)q^m+v(m)$ for some periodic functions $u$ and $v$
of period $d$.

Suppose that $g(x)$ is primitive, so we can take $d=q^d-1$. Let $g(x)$
and $h(x)$ be polynomials of degree less than $d$. If
$$\frac{g(x)}{g(x)}=\frac{b_0+b_1x+\cdots+b_{d-1}x^{d-1}}{1-x^d},$$
then for some $0\leq j\leq d-1$ we have
$$\frac{h(x)}{g(x)}=\frac{b_j+b_{j+1}x+\cdots+b_{d-1}x^{d-j-1}+b_0x^{d-j}
  +\cdots+b_{j-1}x^{d-1}}{1-x^d}.$$
Moreover, all elements of $\ff_q^{*}$ occur equally often among
$b_0,b_1,\dots,b_{d-1}$, while $0$ occurs one fewer times.  Hence each
$h_{im}$ has $q^{d-1}$ coefficients equal to $\alpha\in\ff_q^{*}$
(and $q^{d-1}-1$ coefficients equal to $0$).  Thus the number of
coefficients of $h_{im}$ equal to $\alpha$ has the form
$q^{m+d}/d=q^{m+d}/(q^d-1)$ plus a periodic term. Summing over $0\leq
i\leq d-1$ gives
$$N_{\alpha}(m)=\frac{dq^{m+d}}{q^d-1}+\text{periodic term,}$$
and the proof follows. \qed

\begin{example}
 \be
\item[(a)] Write $[a_0,a_1,\dots,a_{k-1}]$ for the periodic function
$p(m)$ on $\zz$ satisfying $p(m)=a_i$ for $m\equiv i\modd{k}$.
Let $q=2$, $g(x)=1+x+x^2+x^3+x^4$, and $c=1$. The polynomial $g(x)$ is
irreducible over $\zz_2$ but not primitive. It can then be
computed that
 \bea N_1(m) & = &
 2^{m+1}-\frac25(-2)^m+\frac15[-3,1,3,-1].\nonumber\\
   & = & \frac 15[8,12]2^m+\frac 15[-3,1,3,-1]. \label{eq2} \eea
In fact, in Theorem~\ref{thm:1} we can take $d=4$. For $m$ even, the
Taylor series expansion (at $x=0$) of each $h_{im}(x)/g(x)$ has two
coefficients within each period equal to $1$. Hence in this case
$$N_1(m)=4\cdot\frac25\cdot 2^m+\cdots=\frac852^m+\cdots .$$
If $m$ is odd, then the expansions of $h_{0m}(x)/g(x)$ and
$h_{3m}(x)/g(x)$ have two coefficients within each period equal to
$1$, while $h_{1m}(x)/g(x)$ and $h_{2m}(x)/g(x)$ have four such
coefficients equal to $1$. Hence in this case
 $$ N_1(m)=\left(\frac25+\frac45+\frac45+\frac25\right)2^m+\cdots=
    \frac{12}52^m+\cdots, $$
agreeing with equation \eqref{eq2}.
 \item[(b)] Let $q=2$, $g(x)=1+x^2+x^5, c=1$, and $\alpha=1$. Then
$g(x)$ is primitive, and we have
$$N_1(m)=\frac{80}{31}2^m+\frac1{31}[-49,-67,-41,11,-9].$$
\item[(c)] Let $q=2$, $g(x)=1+x+x^3+x^4+x^5, c=1$, and
$\alpha=1$. Then $g(x)$ is primitive, and we have
$$N_1(m)=\frac{80}{31}2^m+\frac1{31}[-49,-5,-41,11,-9].$$
Note that $u(m)=80/31$ for both (b) and (c), as guaranteed by part (d)
of the theorem, but the periodic terms $v(m)$ differ
(though only for $n\equiv1\modd{5}$).
\item[(d)] To illustrate that equation~\eqref{eq1} need not hold for
all $m\geq 0$, let $q=2, g(x)=(1+x^2+x^5)^3$, and $\alpha=1$. Then
$$N_1(m)=\begin{cases} 1,\qquad m=0\\9,\qquad m=1\\
\frac{168}{31}2^m+\frac1{31}[297,-243,-393,-507,-177],\qquad m\geq
2.\end{cases}$$
\item[(e)] Some examples for $q=2$ and $c=3$: first let
$g(x)=1+x+x^2+x^3+x^4$. Then
$$N_1(m)=2^{m+1}-\frac14(-2)^m+\frac15[11,3,-11,-3],\qquad m\geq 2.$$
If $g(x)=1+x^2+x^5$, then
$$N_1(m)=\frac{60}{31}2^m+\frac1{31}[33,-27,-147,-201,-123],\qquad m\geq 2.$$
If $g(x)=1+x+x^2+x^3+x^4+x^5$, then
$$N_1(m)=\frac{60}{31}2^m+\frac1{31}[-153,35,-85,-77,-61],\qquad m\geq 2.$$
\item[(f)] Two examples for $q=3$ and $c=1$. Let $g(x)=2+x+x^2$, a
primitive polynomial. Then for $m\geq0$,
\begin{align}
&N_1(m)=\frac343^m+\frac12-\frac14(-1)^m\\
&N_2(m)=\frac343^m-\frac12-\frac14(-1)^m.\end{align}
Let $g(x)=2+x^2+x^3$, an irreducible but not primitive
polynomial. Then for $m\geq 0$,
\begin{align}
&N_1(m)=\frac{18}{13}3^m+\frac1{13}[-5,11,7]\\
&N_2(m)=\frac9{13}3^m-\frac1{13}[9,14,3].\end{align}
\item[(g)] A class of function for which $u(m)$ is independent of
$m$. Let $q=2, c=1, \alpha=1$ and $g(x)=1+x^{k-1}+x^k$.

If $k=2^h$, then
$$u(n)=\frac{k(3^h-1)}{k^2-1}.$$
If $k=2^h+1$, then
$$u(n)=\frac{k(k-2)(3^h+1)}{2^{3h}-1}.$$
If $k=2^h-1$ and  $d_1,\dots,d_r$ are the degrees of the irreducible
factors of $1+x^{h-1}+x^h\in\zz_2[x]$, and we set
$\delta_h=\mathrm{lcm}\{2^{d_1}-1,\dots,2^{d_r}-1\}$, then
$$u(n)=\frac{k2^{\delta_h-1}}{2^{\delta_h}-1}.$$
 \ee
 \end{example}

Let us consider some examples involving arbitrary powers $g(x)^n$ of
the polynomial $g(x)$.

Let $N_{\alpha}(g)$ denote the number of coefficients of $g(x)$ that
are equal to $\alpha$ and $N(g)$ the total number of nonzero
coefficients of $g(x)$.  Note that $N(g)+N_0(g)=1+\deg g$. Write the
$p$-ary digits of $m\in\nn$ as $\langle m_0,m_1,\dots,m_s\rangle$
so that $m=m_0+m_1p+\cdots+m_sp^s$.

\begin{proposition} \label{prop2.3} {\rm{(a)}} Let $g(x)=(1+x+\cdots
  + x^{p-1})^n\in\zz_p[x]$ for a prime $p$ and $n\in\pp$. Then
  the coefficient of $x^k$ in $g(x)$ is given by
$$[x^k]g(x)=(-1)^k\binom{pn-n}k, \qquad 0\leq k\leq pn-n.$$
{\rm{(b)}} If $\,\langle b_0,b_1,\dots,b_m\rangle$ are the $p$-ary
digits of $(p-1)n$,
then the number of coefficients of $g(x)$ not divisible by $p$ is
 \beq N(g)=\prod_{i=0}^m(1+b_i). \label{eq:ngstar} \eeq
\end{proposition}

\proof Use the trick $1+x+\cdots+x^{p-1}=(1-x)^{p-1}$ in
$\zz_p[x]$. Applying the binomial expansion to $(1-x)^{(p-1)n}$
proves (a). On the other hand, Lucas' theorem implies the congruence
 $$\binom{pn-n}{k}\equiv\prod_{i=0}^{m}\binom{b_i}{k_i},$$
where $\langle k_0,k_1,\dots,k_m\rangle$ are the $p$-ary digits of
$k$. A simple counting argument leads to (b). \qed

\begin{example} \label{ex2.4}

(a) Let $p=3$, $g(x)=(1+x+x^2)^n$ and $\langle b_0,\dots,b_m\rangle$
be the ternary digits of $2n$. Then by equation~\eqref{eq:ngstar},
$$N(g)=N_1(g)+N_2(g)=\prod_{i=0}^m(1+b_i).$$

Observe that $(-1)^{k_i}\binom{2}{k_i}=1$ in $\ff_3$ for any 
$k_i=0,1,2$, but 
$(-1)^{k_i}\binom{1}{k_i}=1,2$ or $0$ depending on $k_i=0,1$ or $2$.
In the latter case both $1$ and $2$ are equally attainable.
Obviously $N(g)+N_0(g)=1+\deg g=2n+1$. It follows that $N_1(g)=3^{m+1}, N_2(g)=0,
N_0(g)=2n+1-3^{m+1}$ when no $b_i=1$. In the case that there is some 
$b_i=1$ then $N_1(g)=N_2(g)=\frac12\prod_{i=0}^m(1+b_i)$ (and hence $N_0(g)=2n+1-2N_1(g)$).
\end{example}

%(b) [\textbf{Tewodros: can you indicate the proofs?}] Let $p=5$,
%  $g(x)=(1+x+x^2+x^3+x^4)^n$ and $\langle
%b_0,\dots,b_m\rangle$ the $5$-ary digits of $4n$. We split into cases.

%(i) If some $b_i=3$, then $N_1(g)=N_2(g)=N_3(g)=N_4(g)=\frac14N(g)$.

%(ii) If some $b_i=1$ but no $b_j=3$, then $3N_1=3N_4=6N_2=6N_3=N(g)$.

%(iii) Let $c_j$ denote the number of $j$'s appearing in
%$\langle b_0,\dots,b_m\rangle$ and $\sigma$ be the
%the permutation $(1)(243)$. If some $b_i=2$ but no $b_j=1$ or $3$, then
 %$$ N_{\alpha}(g)=5^{c_4}\sum_{j+1\equiv\sigma(\alpha)\modd{4}}
 % \binom{c_j}{c_j-j}2^j, \qquad\text{for $1\leq\alpha\leq 4$}. $$
%(iv) If some $b_i=4$ and each $b_i\in\{0,4\}$, then $N_1=N(g)$ and
%thus $N_2=N_3=N_4=0$.
%\end{example}

\begin{example} \label{ex2.5} The above proposition does \emph{not}
  apply to the function $g(x)=(1+x+x^2)^n\in\ff_2[x]$.  Still we
  are able to determine the number of odd coefficients of $g(x)$.
  Polynomials of the form $(1+x^{k-1}+x^k)^n\in\ff_2[x]$ should be
  prone to the technique outlined here. Write $\omega(n)=N(g)$.

Suppose $n=2^j(2^k-1)$. Since $(1+x+x^2)^{2^j}=
1+x^{2^j}+x^{2^{j+1}}$ in $\ff_2[\bm{x}]$, we have 
$\omega(n)=\omega(2^k-1)$. Now
  $$ (1+x+x^2)^{2^k-1}=\frac{1+x^{2^k}+
   x^{2^{k+1}}}{1+x+x^2}. $$
It is easy to check that for $k$ odd we have (writing
$i\equiv t (3)$ for $i\equiv t \modd{3}$)
$$(1+x+x^2)^{2^k-1}=\frac{1+x^{2^k}+x^{2^{k+1}}}{1+x+x^2}=
\sum_{\substack{i=0\\ i\equiv 0(3)}}^{2^k-2}
x^i+\sum_{\substack{i=1\\ i\equiv 1(3)}}^{2^{k+1}-3}
x^i+\sum_{\substack{i=2^k\\ i\equiv 2(3)}}^{2^{k+1}-2} x^i.$$
It follows that $\omega(2^k-1)=(2^{k+2}+1)/3$. Similarly, when $k$ is
even we have
$$(1+x+x^2)^{2^k-1}\equiv\frac{1+x^{2^k}+x^{2^{k+1}}}{1+x+x^2}\equiv
\sum_{\substack{i=0\\i\equiv 0(3)}}^{2^{k+1}-2} x^i+
\sum_{\substack{i=1\\i\equiv 1(3)}}^{2^{k}-2}
x^i+\sum_{\substack{i=2^k+1\\i\equiv 2(3)}}^{2^{k+1}-3} x^i\modd{2}.$$
Hence in this case $\omega(2^k-1)=(2^{k+2}-1)/3$. Now any positive
integer $n$ can be written uniquely as $n=\sum_{i=1}^r2^{j_i}(2^{k_i}-1)$,
where $k_i\geq 1, j_1\geq 0$, and $j_{i+1}>j_i+k_i$. We are simply
breaking up the binary expansion of $n$ into the maximal strings of
consecutive $1$'s. The lengths of these strings are $k_1,\dots,k_r$. Thus
$$(1+x+x^2)^n\equiv
\prod_{i=1}^r(1+x^{2^{j_i}}+x^{2^{j_i+1}})^{2^{k_i}-1}\modd{2}.$$
The key observation is: there are \emph{no cancellations} among the
coefficients when we expand this product since $j_{i+1}>j_i+1$. Hence
 $$\omega(n)=\prod_{i=1}^r\omega(2^{k_i}-1)=
 \prod_{i=1}^r\frac{2^{k_i+2}+(-1)^{k_i+1}}3.$$
Take for instance $n=6039$ with binary expansion $1011110010111$. The
maximal strings of consecutive $1$'s have lengths $1,4,1$ and $3$.
Hence $\omega(6039)=\omega(1)\omega(15)\omega(1)\omega(7)=3\cdot
21\cdot 3\cdot 11=2079$.
%Notice that $\omega(15673)=\omega(6039)$, but
%this is not surprising since permuting the block of strings does not
%alter the computation.

\end{example}

\section{Other multivariate polynomials over $\ff_p[\bm{x}]$ }
\label{sec3}
For the remainder of this paper we make no attempt to be systematic,
but rather confine ourselves to some interesting examples.

Theorem~\ref{thm:fxkn} deals with coefficients of $f(\bm{x})^n$ 
over $\ffq$. We discuss certain cases for which we can be more explicit. We also give some examples of counting coefficients over $\ffq$ of
class of polynomials not of the form $f(\bm{x})^n$ for fixed $f(\bm{x})$.

\begin{example} \label{ex3.1}  The Pascal triangle
modulo an integer $d$ (in the present context, the coefficients
of $(1+x)^n \modd{d}$) receives a good discussion and further references
by Allouche and Shallit \cite[Chapter~14, Section 14.6]{allsha}.

\end{example}

\begin{theorem} \label{thm3.2}
Suppose $H_n(\bm{x})=\prod_{i=1}^n(1+x_i+x_{i+1})\in\ff_p(\bm{x})$
with $\bm{x}=(x_1,\dots,x_{n+1})$ and $p$ a prime. Then we have
$$\sum_{n=0}^{\infty}N(H_n)z^n=\frac{1-z^p}{(1-z)^2-z(1-z^p)}.$$
\end{theorem}

\proof The generating function can be simplified to
$(1+z+\dots+z^{p-1})/(1-2z-z^2-\dots -z^p)$. Thus we need to
verify that the sequence $N(H_n)$ satisfies the recurrence relation
$$t_{n+p}-2t_{n+p-1}-t_{n+p-2}-\dots -t_{n+1}-t_n=0,$$ with
$t_k=N(H_k)$. To this end write
$H_{n+p}=H_p\cdot\prod_{i=1}^n(1+x_{i+p}+x_{i+p+1})$. Clearly it is
enough to prove the case $n=0$; that is,
$t_p=2t_{p-1}+t_{p-2}+\dots+t_1+t_0$. Note that $t_0=1$. We have

\begin{align}
H_p&=(1+x_p+x_{p+1})H_{p-1}\nonumber\\
&=x_{p+1}H_{p-1}+(1+x_p)H_{p-1}\nonumber\\
&=x_{p+1}H_{p-1}+(1+x_{p-1}+x_p^2)H_{p-2}+x_p(2+x_{p-1})
H_{p-2}\nonumber\\
&=x_{p+1}H_{p-1}+\widetilde H_{p-1}+x_p(2+x_{p-1})H_{p-2},
\nonumber\end{align}
where $\widetilde H_{p-1}=(1+x_{p-1}+x_p^2)H_{p-2}$ and
$N(x_{p+1}H_{p-1})=N(\widetilde H_{p-1})=N(H_{p-1})$ (just replace
$x_p\rightarrow x_p^2$).

Continuing the above process we arrive at
 \beq H_p=x_{p+1}H_{p-1}+\widetilde H_{p-1}+\widetilde H_{p-2}+\cdots
 +\widetilde
   H_1+\widetilde H_0, \label{eqsss} \eeq
where
 $$ \widetilde H_0:=x_p\cdots
    x_2(p+x_1)H_0\equiv x_p\cdots x_1H_0\modd{p}, $$ and
 $$ \widetilde
 H_{p-k}:=x_p\cdots x_{p-k+2}(k+kx_{p-k}+x_{p-k+1}^2)H_{p-k-1}$$
for $1\leq k<p$.

In $\ff_p(\bm{x})$ the map $\beta\rightarrow
\beta^2/k$ is bijective and hence $N(\widetilde
H_{p-k})=N(H_{p-k})$. Since the terms in equation~\eqref{eqsss} are
mutually exclusive, a straightforward counting completes the
argument. \qed

\medskip
\textsc{Note.} When $p=2$, it follows from \cite[Exam.~4.1.2]{18} that
$N(H_n)$ is the number of self-avoiding lattice paths of length $n$
from $(0,0)$ with steps $(1,0)$, $(-1,0)$, or $(0,1)$.

\begin{corollary} \label{cor3.3}
If $p\geq 3$ and $h_n(x)=(1+x+x^p)^{(p^n-1)/(p-1)}\in\ff_p[x]$,
then
$$\sum_{n=0}^{\infty}N(h_n)z^n=\frac{1-z^p}{(1-z)^2-z(1-z^p)}.$$
\end{corollary}

\proof In Theorem~\ref{thm3.2}, reindex $x_i$ as $x_{i-1}$ so that
$\bm{x}=(x_0,x_1,\dots,x_n)$; put $x_i=x^{p^i}$ and apply
$a^p+b^p=(a+b)^p \modd{p}$.  
The map between the new $H_n(\bm{x})$ and $h_n(x)$ defined by 
$$ x_{i_1}x_{i_2}\cdots x_{i_r}\rightarrow x^{p^{i_1}+p^{i_2}+\cdots 
   +p^{i_r}} $$
is easily checked to be a bijection due to the uniqueness (up to
permutation of the terms) of $p$-ary digits in the field $\ff_p$. The
result follows from Theorem~\ref{thm3.2}.   
\qed

\section{Counting integer coefficients in shifting products}
\label{sec4}

The motivation for the next discussion comes from the simple formula
 \beq N((x_1+x_2+\cdots+x_r)^s)=\binom{r+s-1}s, \label{eq3} \eeq
the number of $s$-combinations of an $r$-element
set allowing repetition.
Denote $\binom{r+s-1}s$ by $\left(\bsp\binom{r}s\bsp\right)$.

Let $\lambda=(\lambda_1,\dots,\lambda_n)$ be an integer partition, so
$\lambda_1\geq \cdots\geq\lambda_n\geq 0$. Set $\lambda_{n+1}=0$. We
will identify $\lambda$ with its \emph{Young diagram} \cite[p.~29]{18}
consisting of $\lambda_i$ left-justified squares in the $i$th row.

\medskip
 \emph{Notation.}  For $m\in\nn$, write
 $\bm{x}^{(m)}=\sum_{j=1}^mx_j$, with $\bm{x}^{(0)}=0$.  If
$\{\alpha_1,\alpha_2,\dots\}$ is a sequence, then the \emph{forward
difference} is $\Delta \alpha_i=\alpha_{i+1}-\alpha_i$.

\medskip
We now find a generalization of equation~\eqref{eq3}.

\begin{lemma} \label{lemma4.2}
\be\item[(a)] The number of distinct monomials in the
product
$\Omega_n(\lambda)=\prod_{i=1}^n\bm{x}^{(\lambda_i)}$ is equal to
$$N(\Omega_n(\lambda))=\sum_{\bms{k}\in
  K_n}\prod_{i=1}^n\left(\bsp\binom{-\Delta\lambda_i}{k_i}\bsp\right),$$
where $K_n=\{\bm{k}=(k_1,\dots,k_n)\in\nn: k_1+\cdots
+k_i\leq i; \sum k_j=n\}$.

\item[(b)] There holds the recurrence
$$N(\Omega_n(\dots,\lambda_i+1,\dots))=N(\Omega_n(\dots,\lambda_i,\dots))+
$$
$$ \quad
N(\Omega_i(\lambda_1-\lambda_i,\dots,\lambda_{i-1}-\lambda_i,1))\cdot
N(\Omega_{n-i}(\lambda_{i+1},\dots)).$$
\ee
\end{lemma}

\proof Divide up the Young diagram into compartments by drawing
vertical lines alongside the edges $\lambda_i$. The vertices then
generate the set $K_n$ while side-lengths take the form
$\lambda_i-\lambda_{i+1}=-\Delta\lambda_i$.  Now apply  
equation~\eqref{eq3} to each of the disjoint rectangular blocks. \qed

\medskip
\textsc{Remark.} The set $K_n$ is easily seen to be in bijection with
Dyck paths of length $2n$ \cite[Cor.~6.2.3(iv)]{ec2} and hence has
cardinality $C_n$, a Catalan number. The elements of $K_n$ are called
\emph{$G$-draconian sequences} by Postnikov \cite{15}.

\medskip
A very special case involves the ubiquitous Catalan number $C_n$.

\begin{corollary} \label{cor4.4}
Let $\lambda=(n,n-1,\dots,1)$. Then

\be\item[(a)] $N(\Omega_n(\lambda))=C_n$.

%\item[(b)] The weakly-increasing indices of each monomial encode the
%$(1,0)$-steps of the Dyck path at the corresponding heights.

\item[(b)] $C_n=\sum_{j\geq 1}(-1)^{j+1}\binom{n+2-j}jC_{n-j}$.
\ee
\end{corollary}

\proof (a) Since $\Delta\lambda_i=-1$ each product weight equals 1,
therefore $N(\Omega_n(\lambda))=N(K_n)=C_n$.

(b) Consider $\Omega_n(\lambda)$ as lattice boxes in the
region $\{(x,y): y\geq x\geq 0\}$, and let $a(n,j)$ be the number of 
lattice paths from $(0,0)$ to $(n,j)$. This results in the system
$$a(n,j)=\begin{cases} a(n,j-1)\qquad \text{if $j=n$}\\
a(n,j-1)+a(n-1,j)\qquad \text{if $1\leq j<n$},\end{cases}$$
with the conditions $a(0,0)=1$ and $a(n,j)=0$ whenever $j>n$ 
or $j<0$. By construction and part (a), we have 
$a(n,n)=a(n,n-1)=C_{n-1}$. Combining with the relation 
$a(n,n-(j+1))=a(n,n-j)-a(n-1,(n-1)-(j-1))$ and induction, 
one can show that
$$a(n,n-j)=\sum_{i\geq 1}(-1)^{i+1}\binom{j+2-i}iC_{n-i}.$$
Using $\sum_{j=0}^{n-1}a(n,j)=C_n$ and simple
manipulations completes the proof. \qed

Next we list some interesting connections between the result in
Lemma~\ref{lemma4.2} and several other enumerations arising in recent
work by different authors.

\subsection{Lattice paths under cyclically shifting boundaries}
\label{subsec4.1.1}

Chapman \emph{et al.}\ \cite{3} and Irving-Rattan \cite{11} have
enumerated lattice paths under \emph{cyclically shifting} boundaries.
For notations and terminology refer to \cite[Thm.~15]{11}.

\begin{theorem} \label{thm4.1.1} {\rm{(Irving-Rattan)}} Let
  $s,t,n\in\pp$.  Let $U$ and $R$ denote up and right
  steps, respectively. Then there are $\frac1n\binom{(s+t)n-2}{n-1}$ 
  lattice paths from $(0,0)$ to $(sn-1,tn-1)$ with steps $U$ and $R$ 
  lying weakly beneath $U^{t-1}(R^sU^t)^{n-1}R^{s-1}$.
\end{theorem}

Setting $s=t$ yields a result of Bonin-Mier-Noy \cite[Thm.~8.3]{2}, as
follows.

\begin{corollary} \label{cor4.1.2} Let $n$ and $t$ be positive
integers. Then there are
$tC_{nt-1}$ lattice paths from $(0,0)$ to $(nt-1,nt-1)$ with steps $U$ and $R$ lying weakly beneath $U^{t-1}(R^tU^t)^{n-1}R^{t-1}$.
\end{corollary}

Applying Lemma~\ref{lemma4.2} it is possible to give a generating
function reformulation of Theorem~\ref{thm4.1.1}.

\begin{corollary} \label{cor4.1.3}
Let $Z_{n,s,t}$ be the polynomial
 $$Z_{n,s,t}=\left(\sum_{i=1}^{sn}x_i\right)^{t-1}
\prod_{j=1}^{n-1}\left(\sum_{i=1}^{sj}x_i\right)^t.$$
Then the set of distinct monomials in $Z_{n,s,t}$ is equinumerous 
with
lattice paths from $(0,0)$ to $(sn-1,tn-1)$ lying weakly beneath
$U^{t-1}(R^sU^t)^{n-1}R^{s-1}$. Moreover,
$$ N(Z_{n,s,t})=\sum_{\bms{k}\in
  L_{n,t}}\prod_{i=1}^n\binom{k_i+s-1}{k_i}=
\frac1n\binom{(s+t)n-2}{n-1},$$
where $L_{n,t}:=\{\bm{k}=(k_1,\dots,k_n)\in\nn:
k_1+\cdots +k_j\leq tj-1; \sum k_j=tn-1\}$.
%then N(G_{n,s,t})=\frac1n\binom{(s+t)n-2}{tn-1}$.
\end{corollary}

\proof
For convenience, write each factor $x_1+\cdots +x_{\mu}$ in decreasing
order $x_{\mu}+\cdots +x_1$. Treat this as a horizontal box of $\mu$ 
unit squares. Then encode the factors $Z_{n,s,t}$ by way of 
stacking up smaller boxes on top of larger ones, so that everything is right-justified. There will be $n$ rectangular blocks with the 
bottom one of size $sn\times(t-1)$ while the topmost has dimensions $s\times t$. Now the bijection with the lattice paths weakly beneath $U^{t-1}(R^sU^t)^{n-1}R^{s-1}$ is most natural and apparent.

The second assertion results from applying Lemma~\ref{lemma4.2} and
Theorem~\ref{thm4.1.1}. There is only a slight alteration
(simplification) done to the left-hand side. Namely, for $1\leq j\leq n-1$ 
let $I_j:=\{jt,\dots,(j+1)t-1\}$ and choose $\lambda$ to be
$$\lambda_i:=\begin{cases} s(n-j),\qquad i\in I_j\\ sn, \qquad 1\leq i\leq
t-1.\end{cases}$$ 
According to Lemma~\ref{lemma4.2} the underlying set will be $K_{nt-1}$ with 
elements denoted by $\bm{k}=(k_1,\dots,k_{nt-1})$.
By direct calculation we find $\Delta\lambda_i=-s$ for
$i=t-1,2t-1,\dots,nt-1$ and $\Delta\lambda_i=0$ otherwise. 
In the latter case, if the corresponding $k_i\neq 0$ then the related
binomial term vanishes. On the other hand, if such $k_i=0$ then the
binomial contribution is 1. Dropping off the elements $\bm{k}$ with
zero outputs (hence redundant) offers a large reduction on the
set $K_{nt-1}$. Hence the relevant set to sum over becomes $L_{n,t}$. \qed

Write $\lambda\vdash n$ to denote that
$\lambda=(\lambda_1,\lambda_2,\dots)$ is a partition of $n\geq 0$. We
also write $\lambda=\langle 1^{m_1},2^{m_2},\dots\rangle$ to denote
that $\lambda$ has $m_i$ parts equal to $i$.
%Let $\pi$ be the \emph{partition} function on the integers. For
%$u\in\pi(n)$, write $u=[b_1^{e_1},\dots,b_j^{e_j}]$ and its
%\emph{part-length} is $e_1+\cdots +e_j$ so that
%$e_1b_1+\cdots+e_jb_j=n$.
Certain specialized values in
Corollary~\ref{cor4.1.3} produce the following identities.

\begin{corollary} \label{cor4.1.4} With notation as in
Corollary~\ref{cor4.1.3}, if $t=1$ then we have
$$\sum_{\lambda=\langle 1^{m_1},2^{m_2},\dots\rangle\vdash n-1}
\binom{n+1}{m_1+\cdots+m_j}\binom{m_1+\cdots+m_j}
{m_1,\cdots,m_j}\prod_{i\geq 1}\binom{i+s-1}{i}^{m_i} $$
$$ \qquad\qquad=\binom{(s+1)n-2}{n-1}.$$ If instead $s=1$, then
$L_{n,t}:=\{(a_1,\dots,a_n)\in\pp^n_{\geq 0}: a_1+\cdots+a_k\leq
tk-1; \sum a_i=tn-1\}$ is of cardinality
$$ \#L_{n,t} =\frac1n\binom{(t+1)n-2}{n-1}.$$
\end{corollary}

 \subsection{The PS-polytope} \label{subsec4.2}
Let $t_1,\dots,t_n\geq 0$. In \cite{14} Pitman and Stanley discuss the  $n$-dimensional polytope
  $$ \hspace{-2in}\Pi_n(t_1,\dots,t_n) = $$ \vspace{-.3in}
  $$ \{\bm{y}\in\rr^n\st \text{$y_i\geq 0$
and $y_1+\cdots +y_i\leq t_n+\cdots+t_{n-i+1}$ for all $1\leq i\leq
n$}\}.$$
We call this a \emph{PS-polytope.} One of the results in \cite{14}
concerns the total number $\#\Pi_n$ of lattice points in $\Pi_n$ when
each $t_i\in\nn$:
\beq
\#\Pi_n(t_1,\dots,t_n)=\sum_{\bms{k}\in
K_n}\left(\bsp\binom{t_n+1}{k_n}\bsp\right)\prod_{i=1}^{n-1}
\left(\bsp\binom{t_i}{k_i}\bsp\right),\label{pstan} \eeq
where $K_n$ is as in Lemma~\ref{lemma4.2}.

It turns out that the enumeration in Lemma~\ref{lemma4.2} coincides
with that of \cite{14}. The next result makes this assertion
precise; the proof is immediate.

\begin{corollary} \label{cor4.2.1} Let $\bm{t}=(t_1,
\dots,t_n)\in\pp^n$ and $\lambda_i=t_i+t_{i+1}+\cdots+t_n$. Then the
number of monomials in
  $\Omega_n(\lambda):=\prod_{i=1}^n\sum_{j=1}^{\lambda_i}x_j$ equals
  the number of lattice points $\#\Pi_n(t_1,t_2,\dots,t_n-1)$ in the
  PS-polytope $\Pi(t_1,\dots,t_{n-1},t_n-1)$.
\end{corollary}

The term $t_n-1$ can be symmetrized as follows. Consider the so-called
\emph{trimmed generalized permutohedron,} introduced by Postnikov
\cite{15} as the Minkowski sum
$$P_G^{-}(\bm{t})=t_1\Delta_{[n+1]}+t_2\Delta_{[n]}+\cdots+t_n\Delta_{[2]}$$
of the standard \emph{simplices} $\Delta_{[i]}=$conv$(e_1,\dots,e_i)$, where $[i]=\{1,\dots,i\}$, the $e_k$'s are the coordinate vectors in $\rr^i$, and conv denotes convex hull. With these objects defined, 
then the two polytopes are related as
$$P_G^{-}(t_1,\dots,t_n)+\Delta_{[n+1]}=\Pi_n(t_1,\dots,t_n).$$
Hence $N(P_G^{-}(\bm{t}))$ is exactly the count on the distinct monomials
of $\Omega_n(\lambda)$. 

When $t_1=t_2=\dots=t_n+1=t$, a direct proof of Pitman-Stanley's result
~\eqref{pstan} can be given.

\begin{corollary} \label{cor4.2.2} Suppose $\bm{t}=(t,t,\dots,t-1)$.
  Then
$$\#\Pi_n(\bm{t})=\sum_{\bms{k}\in
  K_n}\prod_{i=1}^n\binom{k_i+t-1}{k_i} = \frac 1n
  \binom{(t+1)n-2}{n-1}, $$
the number of lattice paths from $(0,0)$ to $(n-1,nt-1)$ with steps $U$ and $R$ lying weakly beneath $U^{t-1}(RU^t)^{n-1}$ (or 
equivalently $(R^tU)^{n-1}R^{t-1})$.
\end{corollary}

\proof Combining Lemma~\ref{lemma4.2} together with
Corollary~\ref{cor4.1.3} we write
 $$\sum_{\bms{k}\in L_{n,t}}\prod_{i=1}^n\binom{k_i+s-1}{k_i}=
  \frac1n\binom{(s+t)n-2}{n-1}, $$
where $L_{n,t}=\{\bm{k}=(k_1,\dots,k_n)\in\nn^n\vert
k_1+\dots +k_j\leq tj-1; {} \sum k_i=tn-1\}.$

Taking note of the symmetry in $s$ and $t$ (evident from the right
side), compute first at $(t,s)=(1,t)$ and then at
$(t,s)=(t,1)$. The outcome is
 \beq \sum_{L_{n,1}}\prod_{i=1}^n\binom{k_i+t-1}{k_i}=
\sum_{L_{n,t}}\prod_{i=1}^n\binom{k_i+1-1}{k_i}=\frac1n\binom{(t+1)n-2}{n-1}.
 \label{eqstar} \eeq
Observe that the middle term in equation~\eqref{eqstar} is simply the
cardinality of $L_{n,t}$ and by Rado's result \cite{16}
on permutohedrons \cite[Proposition~2.5]{15}, we have 
$L_{n,t}=\Pi_n(\bm{t})$. Also note that
$\#L_{n,1}=\#K_{n-1}=C_{n-1}$. The proof follows. \qed

\subsection{Noncrossing matchings} \label{subsec4.3}

The discussion below has its roots in the following scenario: if
$2n$ people are seated around a circular table, in how many ways can
all of them be simultaneously shaking hands with another person at the
table in such a way that none of the arms cross each other? Answer:
the Catalan number $C_n$.

It is convenient to formulate the above question for a circularly
arranged points $OXOX\cdots OX=(OX)^n$ of $O$'s and $X$'s where only
opposite symbols can be connected. Call the desired goal
\emph{noncrossing matchings}.  Mahlburg-Rattan-Smyth-Kemp \cite{12} have
  extended the concept in a more general setting, i.e., for a string of
  the type $O^{m_1}X^{m_1}\cdots O^{m_n}X^{m_n}$. The next statement
  captures a seemingly nonobvious coincidence.

\begin{corollary} \label{cor4.3.1} Let $0\leq m_1\leq\cdots\leq m_n$
  be integers and $\lambda=(\lambda_1,\dots,\lambda_n)$ where
  $\lambda_i=m_n+\dots+m_i$.  Then the number of noncrossing matchings
   for a string of type $O^{m_1}X^{m_1}\cdots O^{m_n-1}X^{m_n-1}$
   equals the number of monomials in $\Omega_n(\lambda)$.
   Equivalently,

\beq \sum_{\bms{k}\in K_n}\binom{m_n}{k_n}\prod_{i=1}^{n-1}
\binom{m_i+1}{k_i}=\sum_{\bms{k}\in K_n}\prod_{i=1}^n 
\binom{m_i+k_i-1}{k_i}. \label{eq:mrsk} \eeq
\end{corollary}

\proof The left-hand side of the equation above is precisely the
enumeration found in \cite{12} while, the right-hand side is from
Lemma 4.2.  Thus we only need to prove the identity
\eqref{eq:mrsk}.

Given the polytope
$\Pi_n(\bm{m})=\{\bm{y}\in\nn^n: y_1+\cdots +y_i\leq m_n+\cdots
+m_{n+1-i}; 1\leq i\leq n\}$, we know that
$$\#\Pi_n(\bm{m})=\sum_{k\in\bm{k_n}}\binom{m_n+k_n}{k_n}
\prod_{i=1}^{n-1}\binom{m_i+k_i-1}{k_i}.$$
Denoting the interior of $\Pi_n(\bm{m})$ by $\Pi_n^{o}(\bm{m})$ and
applying the \emph{Reciprocity Law} for polytopes \cite[Chapter~4]{1}
yields
\begin{align*}
\#\Pi_n^o(\bm{m})&=(-1)^n\sum_{k\in\bms{k}_n}\binom{-m_n+k_n}{k_n}
\prod_{i=1}^{n-1}\binom{-m_i+k_i-1}{k_i}\\
&=(-1)^n(-1)^{\sum k_i}\sum_{k\in\bms{k}_n}\binom{m_n-1}{k_n}
\prod_{i=1}^{n-1}\binom{m_i}{k_i}\\
&=\sum_{k\in\bms{k}_n}\binom{m_n-1}{k_n} \prod_{i=1}^{n-1}\binom{m_i}{k_i};
\end{align*}
where we used the binomial identity
$\binom{-a+b}{b}=(-1)^b\binom{a-1}b$ and $\sum_{i=1}^nk_i=n$.

On the other hand, $\Pi_n^o(\bm{m})$ is characterized by the
conditions $y_i>0$ and $y_1+\cdots+y_i<m_n+\cdots+m_{n+1-i}$;
equivalently, $y_i^{\prime}:=y_i-1\geq 0$ and
$y_1^{\prime}+\cdots+y_i^{\prime}<(m_n-1)+\cdots+(m_{n+1-i}-1)$ or
$$y_i^{\prime}:=y_i-1\geq 0, \qquad y_1^{\prime}+\cdots+y_i^{\prime}
\leq(m_n-2)+(m_{n-1}-1)+\cdots+(m_{n+1-i}-1).$$
Therefore
\beas \#\Pi_n^o(\bm{m}) & = & \#\Pi_n(m_n-2,m_{n-1}-1,\dots,m_1-1)\\ & = &
\sum_{k\in\bms{k_n}}\binom{m_n+k_n-2}{k_n}
\prod_{i=1}^{n-1}\binom{m_i+k_i-2}{k_i}.
\eeas
Now substituting $m_i\rightarrow m_i+1$ and equating the last two
formulas justifies the assertion. \qed

%the right-hand side of the above equation, which is $N(\Omega_n(\lambda))$.
%So we focus on the left-hand side (LHS). To show this,
%proceed as in Lemma 4.2 but this time the
%vertical stripping be made along the lines $m_i-i+2$ for $1\leq i\leq n-1$,
%and along $m_i-i+1$
%when $i=n$. If $m_i-i+2<0$ or $m_n-n+1<0$ then the contribution is treated as
%zero. On the other hand, the value (LHS)
%is precisely the enumeration found in \cite{7}. \qed

\textsc{Remark.} Although the corollary above has
  been stated for weakly increasing sequences $(m_1,\dots,m_n)$, in
  fact the assertion is valid as a polynomial identity for the
  $n$-tuple $\bm{m}$ of indeterminates. In such generality, however,
  the interpretation in terms of noncrossing matchings will be lost.

\begin{example} \label{ex4.3.3} We consider particular shapes
  $\lambda$ for which we find relations with some known enumerative results.

(a) For $n\geq m$, define the multivariate polynomial
$\phi_{n,m}=\prod_{j=1}^{n-1}\sum_{i=0}^{j+m}x_i$. Then
 $$ N(\phi_{n,m})=\frac{m+2}{2n+m}\binom{2n+m}{n+m+1}, $$
the number of standard Young tableaux of shape $(n+m,n-1)$ as well as
the number of lattice points $\#\Pi_{n-1}(0,1,1,\dots,1,m+2)$.

\medskip
(b) Let
$T_{n,k}=\prod_{j=1}^{n-1}\left(\sum_{i=0}^jx_i\right)^k$ for
$n,k\in\pp$. Then $N(T_{n,k})=\frac1{kn+1}\binom{(k+1)n}n$, a
\emph{Fuss-Catalan number} \cite{Graham}. See also Sloane
\cite{sloane} for a host of other combinatorial interpretations.

\medskip
(c)
Let $M$ and $N$ be two infinite triangular
matrices with entries given for $0\leq j\leq i$ by
\beas M(i,j) & = &
\binom{\binom{i}2-\binom{j}2+\binom{i}1-\binom{j}1+3}{i-j}\\
N(i,j) & = &
\binom{\binom{i}2-\binom{j}2+\binom{i}1-\binom{j}1+2}{i-j}. \eeas
Also let $R:=MN^{-1}$, and denote its elements by
$R(n,k)$. Define the polynomial
$$S_{n,k}=\prod_{j=1}^{n-1}\left(\sum_{i=0}^jx_i\right)^{j+k}.$$
Suppose $\bm{t}=(t_1-1,t_2,\dots,t_{n-1})$ with $t_i=k+n-i$. Then
$N(S_{n,k})=\sum_{\bms{y}\in\Pi_{n-1}(\bms{t})}(1+y_{n-1})=R(n,k)$.
\end{example}

\section{Traveling polynomials in $\ff[\bm{x}]$ } \label{sec5}

Continuing in the same spirit as in the previous section we consider
further multinomials with shifting terms, given by
  \beq W_{j,k,n} = \prod_{i=1}^n (x_{(i-1)j+1}+x_{(i-1)j+2}+\cdots+
          x_{(i-1)j+k}). \label{eq:trapol} \eeq
We call $W_{j,k,n}$ a \emph{traveling polynomial}.  

\begin{theorem} \label{thm5.1}
For fixed $j,k\geq 1$ we have
  \beq \sum_{n\geq 0}N(W_{j,k,n})z^n = \frac{1}{\sum_{h\geq 0}
         (-1)^h \binom{k-j(h-1)}{h}z^h}. \label{eq:thmtr} \eeq
\end{theorem}

\proof
The proof is by the Principle of Inclusion-Exclusion. First note that
if $x_{m_1}x_{m_2}\cdots x_{m_r}$ is a monomial appearing in the
expansion of $W_{j,k,n}$, where $m_1\leq m_2\leq\cdots\leq m_r$, then
we can obtain this monomial by choosing $x_{m_s}$ from the $s$th
factor (i.e., the factor indexed by $i=s$) from the right-hand side of
equation~\eqref{eq:trapol}. Suppose $x_{m_s}=x_{(s-1)j+b_s}$. Then
$b_1 b_2\cdots b_n$ is a sequence satisfying
  \beq 1\leq b_s\leq k,\ \ b_i\geq b_{i-1}-j, \label{eq:bdef} \eeq
and conversely any sequence satisfying equation~\eqref{eq:bdef} yields
a different monomial in the expansion of $W_{j,k,n}$. Hence we want to
count the number $f(n)$ of sequences \eqref{eq:bdef}. Let us call such
a sequence a \emph{valid $n$-sequence}.

For $h\geq 0$ let ${\cal A}_h$ be the set of all sequences
$b_1 b_2\cdots b_n$ satisfying the conditions:
 \begin{itemize}\item $1\leq b_i\leq k$ for $1\leq i\leq n$
  \item $b_i\geq b_{i-1}-j$ for $2\leq i\leq n-h$ 
  \item $b_i<b_{i-1}-j$ for $n-h+2\leq i\leq n$.
  \end{itemize}
In particular, ${\cal A}_0$ is the set of valid $n$-sequences,
so $f(n)=\#{\cal A}_0$, and ${\cal A}_1$ consists of valid
$(n-1)$-sequences followed by any number $1\leq b_n\leq k$. 
For a sequence $\beta=b_1 b_2\cdots b_n$ let 
   $$ \chi(\beta\in{\cal A}_h) =\left\{ \begin{array}{rl}
     1, & \beta\in {\cal A}_h\\ 0, & \beta\not\in {\cal A}_h.
    \end{array} \right. $$
It is easy to see that
  $$ \sum_{h\geq 0}(-1)^h\chi(\beta\in{\cal A}_h) = \left\{
    \begin{array}{rl} 1, & \beta \in{\cal A}_0\\
      0, & \beta \not\in{\cal A}_0. \end{array}\right. $$
It follows that 
  $$ \#{\cal A}_0 = \#{\cal A}_1-\#{\cal A}_2+\#{\cal A}_3-\cdots. $$ 
By elementary and well-known enumerative combinatorics (e.g.,
\cite[Exer.~1.13]{18}) we have 
   $$ \#{\cal A}_h = \binom{k-j(h-1)}{h}f(n-h). $$
Hence 
  $$ f(n) = \sum_{h\geq 1}(-1)^h \binom{k-j(h-1)}{h}f(n-h). $$
This reasoning is valid for all $n\geq 1$ provided we set $f(-1) =
f(-2)= \cdots =f(-(h-1))=0$, from which equation~\eqref{eq:thmtr}
follows immediately.
\qed

\begin{example} \label{ex5.8}
It is rather interesting that even a minor alteration in the traveling
polynomial $W_{j,k,n}$ unveils curious structures. To wit, take the
special case $Q_n:=W_{1,3,n}$. From Theorem~\ref{thm5.1} we gather
that 
\beq 
\sum_{n=0}^{\infty}N(Q_n) z^n=\frac1{1-3z+z^2}, \label{Fib} \eeq
or more explicitly $N(Q_n)=F(2n+2)$, with $F(m)$ the Fibonacci sequence.

Let us now introduce
\beas G_n & = & \prod_{i=1}^n(x_i+x_{i+2}+x_{i+4})\\ 
B_{n,t} & =& \prod_{i=1}^n(1-x_i+x_{i+t}). \eeas
%D_n & = & \prod_{i=1}^n(1-x_i+x_{i+2}). \eeas

(a) We have $N(G_n)=F(n+2)^2-\eta(n)$, where
$\eta(n)=\frac{1-(-1)^n}2$.

\emph{Proof.} Denote
$n^{\prime}=\lfloor{n/2}\rfloor$,
$n^{\prime\prime}=\lfloor{(n+1)/2}\rfloor$.
Writing $G_n$ in the form
$$G_n=\prod_{i\ \text{odd}}(x_i+x_{i+2}+x_{i+4})\prod_{i\ \text{even}}(x_i+x_{i+2}+x_{i+4})$$ 
suggests that (if necessary rename $y_i=x_{2i}$ and $y_i=x_{2i-1}$, 
respectively)
$N(G_n)=N(Q_{n^{\prime}})\cdot N(Q_{n^{\prime\prime}})$. From 
equation~\eqref{Fib} we have
$N(G_n)=F(n+2)^2$ if $n$ is even; $N(G_n)=F(n+1)F(n+3)$
if $n$ is odd. Now invoke \emph{Cassini's formula}
$F(m+1)F(m+3)=F(m+2)^2+(-1)^m$. We also obtain 
 $$ \sum_{n\geq 0}N(G_n)z^n=\frac1{(1-z^2)(1-3z+z^2)}. $$
%Consequently $(\bm{E}^2-1)(\bm{E}^2-3\bm{E}+1)N(D_n)=0$ where $\bm{E}$
%is the unit step \emph{forward shift} operator. After
%some algebraic rearrangement we obtain the identity
%$F(n+1)^2-3F(n)^2+F(n-1)^2=2(-1)^n$.

(b) Recall the fact that the polynomial
  $\prod_{1\leq i<j\leq n}
(x_i-x_j)$ consists of $n!$ nonzero coefficients half of which are
$+1$'s, the other half $-1$'s. In fact, this gives the easiest way to
prove that the number of odd coefficients of $\prod_{1\leq i<j\leq
  n}(x_i+x_j)$  is also $n!$.

A similar phenomenon occurs with $B_{n,t}$: the nonzero coefficients of 
$B_{n,t}$ are all $\pm 1$ with one extra $+1$ than $-1$.
We also have 
\beas \sum_{n\geq 0}N(B_{n,1})z^n & = & \frac{1+z}{1-2z-z^2}\\
\sum_{n\geq
  0}N(B_{n,2})z^n & = &
\frac1{1-z}\left\{\frac{z^2}{1+z^2}+\frac1{1-2z-z^2}\right\}. \eeas
\proof We will only consider $B_{n,1}$; the general case is 
analogous. The
assertions are trivial when $n=0, 1$. Proceed by induction. Start with
the simplification
$$(1-x_n+x_{n+1})(1-x_{n+1}+x_{n+2})=x_{n+1}(1-x_n+x_{n+1})
+(1-x_n+x_nx_{n+1})-x_{n+1}^2.$$
It follows that 
\beq 
B_{n+1,1}=x_{n+2}B_{n,1}+(1-x_n+x_nx_{n+1})B_{n,1}-x_{n+1}^2B_{n-1,1}. 
\label{plusminus} \eeq
Now, the three terms on the right-hand are independent and
$$N((1-x_n+x_nx_{n+1})B_{n-1,1})=N(B_{n,1}).$$
As a consequence $N(B_{n+1,1})=2N(B_{n,1})+N(B_{n-1,1})$ and thereby we 
validate the generating function.

By induction assumption, the polynomials $x_{n+2}B_{n,1}$ and 
$(1-x_n+x_nx_{n+1})B_{n,1}$ each contain
one extra $+1$ than a $-1$; while $-x_{n+1}^2B_{n-1,1}$ has one extra 
$-1$. By equation~\eqref{plusminus} the proof is complete. \qed

\end{example}

We include the next result (Theorem~\ref{thm5.13}) as a generalization of 
the case $j=1$ of 
Theorem~\ref{thm5.1} and as an indication of alternative method of proof towards possible extensions of our work on traveling polynomials. In particular, 
it should be possible to find common generalization of 
Theorem~\ref{thm5.1} and Theorem~\ref{thm5.13}, though 
we have not pursued this question.

Define
 $$ V_{n,k,m}=\prod_{i=1}^n(x_i+x_{i+1}+\cdots +x_{i+k})^m, \qquad
     n,k,m\in\nn. $$
We need a preliminary result before stating Theorem~\ref{thm5.13}.

\begin{lemma} \label{lemma5.12}
Let $k,m$ be nonnegative integers and $A^{(k,m)}$ be the
$(k+1)\times(k+1)$ matrix
given by
$$A_{i,j}^{(k,m)}=\begin{cases} \binom{m-1+i}{m-1}\qquad \text{if
    $j=0$},\\[.2em]
     \binom{m+i-j}{m-1}\qquad \text{if $j\neq 0$};\end{cases} \qquad
   0\leq i,j\leq k.$$
Then $A^{(k,m)}$ has the characteristic polynomial
$$\Theta_{k,m}(\rho)=\sum_{\tau\geq
  0}(-1)^{\tau}\binom{1+(k+1-\tau)m}{\tau}\rho^{k+1-\tau}.$$
\end{lemma}

\proof Fix $m$ and induct on $k$. We want to compute
$\Theta_{k,m}(\rho)=\det(A^{(k,m)}-\rho I)$.  For brevity write
$\Theta_k$ instead of $\Theta_{k,m}$. Observe that apart from the first
column, the matrix $A^{(k,m)}$ is almost Toeplitz with one
super-diagonal. Taking advantage of this fact, extract the determinant
using a repeated application of Laplace expansion along the last
column. We determine the recurrence
 $$\Theta_k=(m-\rho)\Theta_{k-1}+(-1)^k\binom{m+k-1}{m-1}+
 \sum_{i=2}^k(-1)^{i+1}\binom{m-1+i}{m-1}\Theta_{k-i}.$$
Absorb $\Theta_k$ and $m\Theta_{k-1}$ into the sum:
$$0=-z\Theta_{k-1}+(-1)^k\binom{m-1+k}{m-1}+\sum_{i=0}^k(-1)^{i+1}
\binom{m-1+i}{m-1}\Theta_{k-i}.$$
The statement of the lemma is then equivalent to
\begin{align}
0&=\sum_{\tau=0}^k(-1)^{k+1-\tau}\binom{1+(k-\tau)m}{\tau}
\rho^{k+1-\tau}+(-1)^k\binom{m-1+k}{m-1} \nonumber\\
&\ +(-1)^k\binom{m+k}{m-1} \nonumber\\
&\ +\sum_{i=0}^{k+1}(-1)^{i+1}\binom{m-1+i}{m-1}\sum_{j=0}^{k+1-i}
(-1)^{k+1-i-j}\binom{1+(k+1-i-j)m}j\rho^{k+1-i-j}.
\nonumber\end{align}
Replace $i+j\rightarrow \tau$ and rearrange the double summation
\begin{align}
0&=\sum_{\tau=0}^k(-1)^{k+1-\tau}\rho^{k+1-\tau}
\binom{1+(k-\tau)m}{\tau} \nonumber\\
&\ +(-1)^k\binom{m-1+k}{m-1}+(-1)^k \binom{m+k}{m-1}\nonumber \\
&\ +\sum_{\tau=0}^{k+1}(-1)^{k+1-\tau}\rho^{k+1-\tau}
\sum_{j=0}^{\tau}(-1)^{\tau-j+1}\binom{m-1+\tau-j}{m-1}
\binom{1+(k+1-\tau)m}j. \nonumber \end{align}
After comparing coefficients in these polynomials, the problem
is tantamount to proving
 $$\sum_{j=0}^{\tau}(-1)^{\tau-j}\binom{m-1+\tau-j}{m-1}
 \binom{1+(k+1-\tau)m}j=\binom{1+(k-\tau)m}{\tau}.$$
But this formula is a consequence of the Vandermonde-Chu identity
 $$ \sum_{j=0}^{\tau}(-1)^{\tau-j}\binom{m-1+\tau-j}{m-1}
 \binom{x+m}j=\binom{x}{\tau}. $$
\qed

\begin{theorem}\label{thm5.13}
% MASTER THEOREM 5.13.
With $A$ as in Lemma~\ref{lemma5.12} and $\Phi_{\xi}$ being
 the sum of the first column of $A^{\xi}$, we have
 $$\sum_{n\geq 0}N(V_{n,k,m}) z^n=
\frac{\sum_{\nu=0}^{k-1}z^{\nu}\sum_{i=0}^{\nu}(-1)^i
  \binom{1+(k+1-i)m}i\Phi_{\nu-i}}
{\sum_{i\geq 0}(-1)^i\binom{1+(k+1-i)m}iz^i}.$$
In fact, $\Phi_{\xi}=N(V_{\xi,k,m})$ for $0\leq\xi\leq k-1$.
\end{theorem}

 \proof   Given $k, m$ construct the list
 $\{a^{(k,m)}(n,0),\dots,a^{(k,m)}(n,k)\}$ the same way as in
the proof of Corollary~\ref{cor4.4}. Suppress 
$k,m$ and
define the column vector $\Gamma_n=[a(n,0),\dots,a(n,k)]^T$.
Reintroduce the \emph{connectivity matrix} $A^{(k,m)}$, or simply
$A=(A_{i,j})$, from Lemma~\ref{lemma5.12},
%$$A_{i,j}:=\begin{cases} \binom{m-1+i}{m-1}\qquad \text{if $j=0$},\\
%\binom{m+i-j}{m-1}\qquad \text{if $j\neq 0$};\end{cases}
%\qquad
%0\leq i,j\leq k.$$
and apply this result to obtain
$$\sum_{\tau\geq 0}(-1)^{\tau}\binom{1+(k+1-\tau)m}{\tau}A^{k+1-\tau}
=\bm{0}.$$
On the other hand,
\begin{align}
a(n,i)&=\sum_{j=0}^i\binom{m-2+i-j}{m-2}\sum_{c=0}^{j+1}a(n-1,c)
\nonumber\\
&=\sum_{c=0}^{i+1}a(n-1,c)\sum_{j=c-1}^i\binom{m-2+i-j}{m-2}
\nonumber\\
&=\binom{m-1+i}{m-1}a(n-1,0)+\sum_{c=1}^{i+1}\binom{m+i-c}{m-1}a(n-1,c),
\nonumber
\end{align}
where we set $a(\cdot,c)=0$ if $c<0$ or $c>k$. Therefore
$A\Gamma_{n-1}=\Gamma_n$ and hence
 $$\sum_{\tau\geq 0}(-1)^{\tau}\binom{1+(k+1-\tau)m}{\tau}
  \Gamma_{n+k+2-\tau}=\bm{0}.$$
This component-wise sum together with
$N(V_{n,k,m})=\sum_{i=0}^ka(n,i)$ implies that
$$\sum_{\tau\geq 0}(-1)^{\tau}\binom{1+(k+1-\tau)m}
 {\tau}N(V_{n+k+1-\tau,k,m})=0.$$
The denominator of the generating function has been justified and the
proof is complete. \qed

\begin{example} \label{ex5.14}
(a) The case $k=2, m=2$ recovers certain \emph{Kekule numbers} 
$\kappa_n$ \cite[page~302]{4} whose molecular graphs posssess remarkable combinatorial properties often related to the Fibonacci sequence.

(b) Let $k=2$. Recall the \emph{Narayana numbers}
$Y(i,j)=\frac1j\binom{i}{j-1}\binom{i-1}{j-1}$, which enumerate Dyck
paths on $2i$ steps
with exactly $j$ peaks. Then we have
  $$ \sum_{n\geq 0}N(V_{n,2,m})z^n =
  \frac{1+Y(m,2)z}{1-\binom{2m+1}1z+\binom{m+1}2z^2}.$$

(c) If $k=3$ then we have
 $$ \sum_{n\geq 0}N(V_{n,3,m})z^n = \frac{1+\left\{2\binom{m}2+\binom
     {m+1}3\right\}z+Y(m,3)z^2}
    {1-\binom{3m+1}1z+\binom{2m+1}2z^2-\binom{m+1}3z^3}.$$

(d) If $k=4$ then we have
 $$ \sum_{n\geq 0}N(V_{n,3,m})z^n = \frac{1+
  a(m)z+ b(m)z^2+Y(m,4)z^3}
 {1-\binom{4m+1}1z+\binom{3m+1}2z^2-\binom{2m+1}3z^3
 +\binom{m+1}4z^4}, $$
where $a(m)=3\binom{m}2+2\binom{m+1}3+\binom{m+2}4$ and $b(m)=
10\binom{m}3+23\binom{m}4+10\binom{m}5$.

\medskip

(e) In sharp contrast to the above, the one-variable counterpart
(setting $x_i=x^i$) to Theorem~\ref{thm5.13} is a lot simpler.
Namely, if 
$$J_{n,k,m}=\prod_{i=1}^n(1+x^i+x^{i+1}+\cdots+ x^{i+k})^m$$
then
 $$ N(J_{n,k,m})=1+\left\{nk+\binom{n+1}2\right\}m. $$

\end{example}

\begin{example} \label{ex5.15} Let $D_{n,k}$ be the polynomial
$$D_{n,k}(\bm{x},\bm{y})=\prod_{i=1}^n(y_1+\cdots+
y_{i-1}+x_i+\cdots+x_{i+k}).$$  
Notice that
$N(D_{n+1,-1})$ is the polynomial of equation~\eqref{eq:cat}, so its
number of terms is the Catalan number $C_n$. We calculate $N(E_{n,k})$
for $k=0$ and $k=2$.

(i) The case $k=0$. We will be referring to
the discussion and notations from Section~\ref{sec4}. 
Let $[n]=\{1,2,\dots,n\}$ be the integer interval. Given a subset
$A\subset [n]$ associate the multivariate polynomial 
$\Omega(A)= \prod_{a\in A}(y_1+\cdots+y_a)$.

\medskip

Consider the triangular product $\Omega([n])=\prod_{i=1}^n(y_1+\cdots+y_i)$.
By the remark following
Lemma~\ref{lemma4.2}, the monomials in $\Omega([n])$ can be regarded
as Dyck paths of length $2n$. 
From Corollary~\ref{cor4.4} we already know that
$N(\Omega([n])=N(\Omega_n(\lambda))=C_n$, corresponding to
the partition $\lambda=(n,n-1,\dots,1)$. 

\medskip

In the polynomials $D_{n+1,0}(\bm{x},\bm{y})
=\prod_{i=1}^{n+1}(y_1+\cdots+ y_{i-1}+x_i)$
if we replace $x_i\rightarrow y_i$, then clearly $D_{n+1,0}(\bm{x},\bm{y})
=\Omega([n+1])$. Therefore it is natural
to identify the terms in $D_{n+1,0}$ as Dyck paths with ascents (the
edges $x_i$) that are \emph{two-colored}. Such an interpretation and viewing 
the Narayana numbers $Y(n,j)$ as the number of Dyck paths with
exactly $j$ \emph{peaks} \cite{veljan} produces 
the enumeration
\beq 
N(D_{n+1,0})=\sum_{j=0}^n2^jY(n,j).\label{narayana} \eeq
The sum $\sum_{j=1}^n2^j\frac1n\binom{n}j\binom{n}{j-1}$ 
on the right-hand side bears the name \emph{large Schr\"oder
numbers}.

%A similar method as in
%Theorem~\ref{thm5.13} exhibits the number of monomials
%$N(E_{n+1,0})$ being equal 
%With $Y(n,j)$ the Narayana numbers, this amounts to
%$$N(E_{n+1,0})=\sum_{j=0}^n2^jY(n,j).$$
%The sum on the right-hand side bears the name 
%\emph{large Schr\"oder
%numbers} given by
%$$\sum_{j=1}^n2^j\frac1n\binom{n}j\binom{n}{j-1}.$$

%Next, we relate these numbers with the results
%from Section 4 and maintain the notations therein.
%Let $[n]=\{1,2,\dots,n\}$ be the integer interval. Given a subset
%$A\subseteq [n]$, arrange its elements in increasing order as
%$\lambda_{A}=(j_1>j_2>\cdots)$. Assume $j_i=0$ whenever $i>\vert A\vert$.  
%We designate $\vert A\vert$ to be the cardinality of a set $A$.
%Associate the multivariate polynomial 
%$\Omega(A)= \prod_{a\in A}(x_1+\cdots+x_a)$

%Consider the triangular product $\prod_{i=1}^n(x_1+\cdots+x_i)$, which
%is $\Omega_n(\lambda)$ where $\lambda=(n,n-1,\dots,1)$.  Before
%expansion, we view $\Omega_n(\lambda)$ as a product of $n$ rows.
%Suppose now we start crossing out some of these rows. Such an
%elimination can be carried out in $2^{n}$ different ways,
%corresponding to the power set of $[n]$.  

The same analysis actually translates as breaking up the product
$\Omega_n(\lambda)$ into monomials that are cataloged according to how
many $x_i$'s appear in them. We may thus find
 $$N(D_{n+1,0})=\sum_{A\subseteq [n]}N(\Omega(A)).$$
Finally, if we employ Lemma~\ref{lemma4.2}, equation~\eqref{narayana}
and the formula for the Schr\"oder numbers then the outcome
 $$\sum_{A\subseteq[n]}\sum_{\bms{k}\in K_{\#A}}
 \prod_{i=1}^{\#A}\left(\bsp\binom{j_i-j_{i+1}}{k_i}\bsp\right)=
 \sum_{j=1}^n2^j\frac1n\binom{n}j\binom{n}{j-1}$$
is an explicit expression where the sum runs through the power set of
$[n]$.

Convention: when $A$ is the empty set, take the product to be $1$.

\medskip

(ii) The case $k=2$. Interestingly, the values
$\gamma_n=N(D_{n-2,2})$ correspond exactly to the enumeration of
\emph{tandem (unrooted) duplication trees} in gene replication
\cite{8}.

Claim: For $n\geq 3$, we have $\gamma_n=\frac12\nu_n$ where
$\nu_n=\sum_{j\geq 1}(-1)^{j+1}\binom{n+1-2j}j\nu_{n-i}$ with
$\nu_2=1$.

 \proof  Proceed with the first argument in the proof 
of Theorem~\ref{thm5.13}
and put the list $\{a(n,0),\dots,a(n,n-1)\}$, 
so that $\nu_n=\sum_{i=0}^{n-1}a(n,i)$. It is easy to check that this
doubly-indexed sequence $a(n,i)$ satisfies the partial-difference
boundary value problem (BVP)
 $$\begin{cases} a(n,i)=a(n-1,i+1)+a(n,i-1),\\
   a(n,0)=a(n-1,0)+a(n-1,1),\\
   a(n,n-1)=a(n,n-2)=a(n,n-3)=\nu_{n-1}.\end{cases}$$
If $\widetilde\nu$ and $\widetilde{a}$ are the column vectors
\begin{align*}
\widetilde\nu&=[\nu_{n-1},\nu_{n-1},\nu_{n-1},\nu_{n-2},
\nu_{n-3},\dots,\nu_2]^T, \\
\widetilde{a}&=[a(n,n-1),a(n,n-2),\dots,a(n,0)]^T \end{align*}
and $M_{i,j}=(-1)^{i+1}\binom{j+1-2i}{i-1}$ is a matrix, then
the BVP is solved by $M\widetilde\nu=\widetilde{a}$. That means
$$a(n,n-j)=\sum_{i\geq 1}(-1)^{i+1}\binom{j-2i}{i-1}\nu_{n-i}.$$
The claim follows from $\nu_n=\sum_{i=1}^na(n,n-j)$ and
the initial conditions $a(3,0)=a(3,1)=a(3,2)=1$. \qed
\end{example}

\medskip
\textbf{Acknowledgments.}  The authors wish to thank M.
Bousquet-M\'elou, R.  Du, D. E. Knuth, K. Mahlburg, A. Postnikov and C. 
Smyth for valuable discussions.  Richard Stanley's
   contribution is based upon work supported by the National
 Science Foundation under Grant No.~DMS-0604423. Any opinions,
 findings and conclusions or recommendations expressed in this
 material are those of the author and do not necessarily reflect those
 of the National Science Foundation.

\end{document}